\documentclass{article}
\usepackage{amssymb}
\newcommand{\eq}{\begin{equation}}
\newcommand{\en}{\end{equation}}

\newcommand{\giv}{\,|\,}
\newcommand{\boxbar}{\,|\,}
\newcommand{\XX}{{\cal X}}
\newcommand{\FF}{{\cal F}}
\renewcommand{\SS}{S}
\newcommand{\OPi}{{\widetilde{\Pi}}}

\newcommand{\nutil}{\tilde{\nu}}

\newcommand{\re}[1]{\mbox{(\ref{#1})}}
\newcommand{\rem}[1]{\mbox{\rm (\ref{#1})}}

\newcommand{\prob}{\mathbb P}
\newcommand{\ex}{\mathbb E}

\newcommand{\ed}{ \stackrel{d}{=}}

\newcommand{\drift}{{\tt d}}
\def\endpf{\hfill $\Box$ \vskip0.5cm}
\def \proof{\noindent{\it Proof.\ }}

\newtheorem{theorem}{\large Theorem}

\newtheorem{definition}[theorem]{\large Definition}
\newtheorem{corollary}[theorem]{\large Corollary}
\newtheorem{lemma}[theorem]{\large Lemma}

\begin{document}

\title{Regenerative partition structures 
\thanks{
Research supported in part by N.S.F. Grant DMS-0405779
}
}
\author{Alexander Gnedin\thanks{Utrecht University; e-mail gnedin@math.uu.nl}
\hspace{.2cm}
and 
\hspace{.2cm}
Jim Pitman\thanks{University of California, Berkeley; e-mail pitman@stat.Berkeley.EDU} 
\\
\\
\\
\\
}
\date{
\today
\\
}
\maketitle

\centerline{\bf Abstract}
We consider Kingman's partition structures which are regenerative with respect
to a general operation of random deletion of some part.
Prototypes of this class are the Ewens partition structures which 
Kingman characterised by regeneration after deletion of a part chosen by size-biased sampling.
We associate each regenerative partition structure with a corresponding 
regenerative composition structure, which (as we showed in a previous paper) 
can be associated in turn with a regenerative random subset of the positive 
halfline, that is the closed range of a subordinator.
A general regenerative partition structure is thus represented in terms 
of the Laplace exponent of an associated subordinator.
We also analyse deletion properties characteristic of the 
two-parameter family of partition structures. 

\newpage

\section{Introduction and main results}
\label{intro}

This paper is concerned with probability distributions for a random partition 
$\pi_n$ of a positive integer $n$. We may represent $\pi_n$ 
as a sequence of integer-valued random variables
$$
\pi_n = ( \pi_{n,1}, \pi_{n,2}, \ldots ) ~~\mbox{ with }~~ \pi_{n,1} \ge \pi_{n,2} \ge \cdots \ge 0
$$
so $\pi_{n,i}$ is the size of the $i$th largest part of $\pi_n$,
and $\sum_i \pi_{n,i} = n$.
We may also treat $\pi_n$ as a multiset of positive integers with sum $n$,
regarding $\pi_n$ as a random allocation of $n$ unlabeled
balls into some random number of unlabeled boxes, with each box containing
at least one ball.
We call $\pi_n$ {\em regenerative} if it is possible to delete
a single box of balls from $\pi_n$ in such a way that for each
$1 \le x \le n$, given the deleted box contained $x$ balls, the remaining
partition of $n-x$ balls is distributed as if $x$ balls had been deleted
from $\pi_n$ by uniform random sampling without replacement.

\par
To be more precise,
we assume that $\pi_n$ is defined on some probability space 
$(\Omega, \FF, \prob)$ which is rich enough to allow
various further randomisations considered below, including the 
choice of some {\em random part} $X_n \in \pi_n$, meaning that $X_n$
is one of the positive integers in the multiset $\pi_n$ with sum $n$.
The distribution of $\pi_n$ is then specified by some {\em partition probability function}
\eq
\label{pdef}
p(\lambda) := \prob ( \pi_n = \lambda) \qquad ( \lambda \vdash n )
\en
where the notation $\lambda \vdash n$ indicates that $\lambda$ is a partition of $n$. 
The joint distribution of
$\pi_n$ and $X_n$ is determined by the partition probability function $p$ 
and some {\em deletion kernel}
$d = d(\lambda,x), \lambda \vdash n, 1 \le x \le n$, which describes the conditional distribution of $X_n$
given $\pi_n$, according to the formula
\eq
\label{dker}
p(\lambda) d(\lambda, x ) = \prob ( \pi_n = \lambda, X_n = x).
\en
The requirement that $X_n$ is a part of $\pi_n$ makes $d(\lambda,x) = 0$
unless $x$ is a part of $\lambda$, and
\eq
\label{dcon}
\sum_{x \in \lambda} d(\lambda,x ) = 1
\en
for all partitions $\lambda$ of $n$.
Without loss of generality, we suppose further that $\pi_n$ is the 
sequence of ranked sizes of classes of some
random partition $\Pi_n$ of the set $[n]:= \{1, \ldots, n \}$, where 
conditionally given $\pi_n$ all possible values of $\Pi_n$ are equally likely.
Equivalently, $\Pi_n$ is an exchangeable random partition of $[n]$
as defined in \cite{csp}.
For $1 \le m \le n$ let $\Pi_m$ be the restriction of
$\Pi_n$ to $[m]$, and let $\pi_m$ be the sequence of ranked sizes of classes of $\Pi_m$.
We say that the random partition $\pi_m$ of $m$ is {\em derived from $\pi_n$ by random sampling}, 
and call the distributions of the random partitions $\pi_m$ for $1 \le m \le n$ 
{\em sampling consistent}.
A {\em partition structure} is a function $p(\lambda)$ 
as in \re{pdef}
for a sampling consistent sequence of distributions of $\pi_n$ for $n = 1,2, \ldots$. 
This concept was introduced by Kingman \cite{ki78a}, who established a one-to-one
correspondence between partition structures $p$ and distributions for a sequence
of nonnegative random variables $V_1, V_2, \ldots$ with $V_1 \ge V_2 \ge \ldots$ and $\sum_i V_i \le 1$.
In Kingman's {\em paintbox representation} of $p$,
the random partition $\pi_n$ of $n$ 
is constructed as follows 
from $(V_k)$ and a sequence of independent random variables $U_i$ with
uniform distribution on $[0,1]$, where  $(U_i)$ and $(V_k)$  are independent:
$\pi_n$ as in \re{pdef} is defined to be the sequence of ranked sizes of blocks of the partition of $[n]$ 
generated by a random equivalence relation $\sim$ on positive integers, with $i \sim j$ iff either 
$i = j$ or both $U_i$ and $U_j$ fall in $I_k$ for some $k$, where the $I_k$ are some disjoint random 
sub-intervals of $[0,1]$ of lengths $V_k$.
See also \cite{csp} and papers cited there for further background.

\begin{definition}
\label{defreg}
{\em
Call a random partition $\pi_n$ of $n$ {\em regenerative}, 
if it is possible to select a random part $X_n$ of $\pi_n$
in such a way that for each $1 \le x < n$,
conditionally given that $X_n = x$ the remaining partition of $n-x$ is distributed
according to the unconditional distribution of $\pi_{n-x}$ derived from
$\pi_n$ by random sampling.
Then $\pi_n$ may also be called 
{\em regenerative with respect to deletion of $X_n$}, or
{\em regenerative with respect to $d$} if the conditional
law of $X_n$ given $\pi_n$ is specified by a deletion kernel $d$ as in \re{dker}.
Call a partition structure $p$ {\em regenerative} if the corresponding $\pi_n$ 
is regenerative for each $n = 1,2, \ldots$.
}
\end{definition}

According to this definition, $\pi_n$ is regenerative with respect to deletion
of some part $X_n \in \pi_n$ if and only if for each partition $\lambda$ of $n$
and each part $x \in \lambda$,
\eq
\label{def1}
\prob ( \pi _n = \lambda , X_n = x ) = \prob ( X_n = x ) \, \prob ( \pi_{n-x} = \lambda - \{ x \} )
\qquad ( \lambda \vdash n )
\en
where $\lambda - \{ x \}$ is the partition of $n-x$ obtained by deleting the part $x$ from $\lambda$,
and $\pi_{n-x}$ is derived from $\pi_n$ by sampling.
Put another way, $\pi_n$ is regenerative with respect to a deletion kernel $d$ iff
\begin{equation}\label{d-eq}
p(\lambda) d(\lambda,x) = q(n,x)p(\lambda - \{x\})\,,\qquad x\in \lambda
\qquad ( \lambda \vdash n )
\end{equation}
where $p(\mu):= \prob(\pi_m = \mu)$ for $\mu \vdash m$ and $1 \le m \le n$ and
\eq
q(n,x):=\sum_{\{\lambda\vdash n\,:\,x\in \lambda\}} d(\lambda,x)p(\lambda)
=
\prob(X_n = x)
\qquad ( 1\leq x\leq n )
\en
is the unconditional probability that the deletion rule removes a part of size $x$ from $\pi_n$.

\par 
A well known partition structure is obtained by letting
$\pi_n$ be the partition of $n$ generated by the sizes of cycles of a uniformly distributed random permutation
$\sigma_n$ of $[n]$. 
If $X_n$ is the size of the cycle of $\sigma_n$ containing $1$,
then $\pi_n$ is regenerative with respect to deletion of $X_n$, because
given $X_n = x$ the remaining partition of $n-x$ is generated by the cycles of 
a uniform random permutation of a set of size $n-x$.
In this example, the unconditional distribution $q(n,\cdot)$ of $X_n$ is uniform on $[n]$.
The deletion kernel is
\eq
\label{sb}
d(\lambda, x ) = { x a_x (\lambda) \over n}
\qquad (\lambda \vdash n)
\en
where 
$\lambda$ is a partition of $n$ and 
$a_x(\lambda)$ is the number of parts of $\lambda$ of 
size $x$, so $\sum_{x = 1}^n x a_x(\lambda) = n$.
More generally, a part $X_n$ is chosen from a random partition $\pi_n$ of $n$
according to \re{sb} may be called a {\em size-biased part} of $\pi_n$.
According to a well known result of Kingman \cite{ki78a},
if a partition structure is regenerative with respect to deletion of a size-biased part,
then it is governed by the Ewens sampling formula 
\eq
p(\lambda) = {n ! \theta ^\ell \over ( \theta )_{n \uparrow 1 } } \prod_{r} {1 \over r^{a_r} a_r ! }
\en
for some parameter $\theta \ge 0$, 
where $\lambda$ is encoded by its multiplicities $a_r = a_r(\lambda)$ for $r = 1,2, \ldots$, 
with 
\eq
\label{constraints}
\ell=\Sigma\, a_r\,, n=\Sigma \,ra_r\,
\en
and
$$(\theta )_{n\uparrow b}:=\prod_{i=1}^n (\theta +(i-1)b).$$
The case $\theta = 1$ gives the distribution of the partition generated by
cycles of a uniform random permutation.
Pitman \cite{jp.ew,jp.epe} introduced a two-parameter extension of the
Ewens family of partition structures, defined by the sampling formula
\eq
\label{param2}
p(\lambda) = {n ! (\theta )_{\ell \uparrow \alpha } \over ( \theta )_{n \uparrow 1 } } \prod_{r} \left( { ( 1 - \alpha)_{r-1 \uparrow 1 } \over r! } \right) ^ {a_r} { 1 \over a_r ! } 
\en
for suitable parameters $(\alpha,\theta)$, including
\eq
\label{both_pos}
\{ ( \alpha,\theta) : 0\leq\alpha \le 1, \, \theta\geq 0 \}
\en
where boundary cases are defined by continuity.
See \cite{csp} for a review of various applications of this formula. 
The result of \cite[Theorem 8.1 and Corollary 8.2]{gnedinp03} shows that each $(\alpha,\theta)$ partition structure with parameters subject to \re{both_pos}
is regenerative with respect
to the deletion kernel
\begin{equation}\label{dk}
d(\lambda, r)={a_r\over n}\, {(n-r)\tau+r(1-\tau)\over 1-\tau+(\ell-1)\tau}\,,
\end{equation}
where $\tau = \alpha/(\alpha + \theta ) \in [0,1]$, and 
\re{dcon} follows easily from \re{constraints}.
In Section \ref{proofdk} we establish: 

\begin{theorem}
\label{thm2}
For each $\tau\in [0,1]$, 
the only partition structures which are
regenerative with respect to the deletion kernel \rem{dk} are
the $(\alpha,\theta)$ partition structures subject to {\em \re{both_pos}}
with $\alpha/(\alpha+\theta) = \tau$.
\end{theorem}

\par
The following three cases are of special interest:

\paragraph{Size-biased deletion}
This is the case
$\tau=0$: each part $r$ is selected with probability proportional to $r$. 
Here, and in following descriptions, we assume that the parts of a partition
are labeled in some arbitrary way, to distinguish parts of equal size.
In particular, if $\pi_n$ is the partition of $n$ derived from an exchangeable 
random partition $\Pi_n$ of $[n]$, then for each $i \in [n]$ the size $X_n(i)$ 
of the part of $\Pi_n$ containing $i$ defines a size-biased pick from the
parts of $\pi_n$. 
Theorem \ref{thm2} in this case reduces to Kingman's
characterisation of the Ewens family of $(0,\theta)$ partition structures.
Section \ref{related} 
compares Theorem \ref{thm2} with another characterisation of 
$(\alpha,\theta)$ partition structures provided by Pitman \cite{jp.isbp}
in terms of a size-biased random permutation of parts defined by
iterated size-biased deletion.

\paragraph{Unbiased (uniform) deletion}
This is the case $\tau=1/2$: given that $\pi_n$ has $\ell$ parts, each part is
chosen with probability $1/\ell$. Iteration of this operation puts the parts of $\pi_n$ in an exchangeable random order. In this case, the conclusion of Theorem \ref{thm2} is that
the $(\alpha,\alpha)$ partition structures for $0 \le \alpha \le 1$ are the
only partition structures invariant under uniform deletion. 
This conclusion can also be drawn from Theorem 10.1 of \cite{gnedinp03}. As shown 
in \cite{jp.epe,jp.bmpart}, the $(\alpha,\alpha)$ partition 
structures are generated by sampling from the interval partition of $[0,1]$ into
excursion intervals of a Bessel bridge of dimension $2-2 \alpha$. 
The case $\alpha = 1/2$ corresponds to excursions of a standard Brownian bridge.

\paragraph{Cosize-biased deletion}
In the case $\tau=1$, each part of size $r$ is selected 
with probability proportional to the size $n-r$ of the remaining partition.
The conclusion of Theorem \ref{thm2} in this case is that
the $(\alpha,0)$ partition structures for $0 \le \alpha \le 1$ are the
only partition structures invariant under this operation.
As shown in \cite{jp.epe,jp.bmpart}, these partition 
structures are generated by sampling from the interval partition generated by 
excursion intervals of an unconditioned Bessel process of dimension 
$2-2 \alpha$. The case $\alpha = 1/2$ corresponds to excursions of a standard Brownian motion.

\vskip0.5cm
\par
The next theorem, which is proved in Section \ref{frag},
puts Theorem \ref{thm2} in a more general context:

\begin{theorem}
\label{thm1}
For each probability distribution 
$q(n, \cdot\,)$ on $[n]$, there exists a unique joint distribution of
a random partition $\pi_n$ of $n$ and a random part $X_n$ of $\pi_n$ such that 
$X_n$ has distribution $q(n, \cdot)$ and $\pi_n$ is regenerative with respect to deletion of $X_n$.

Let $\pi_m,1 \le m \le n$ be derived from $\pi_n$ by random sampling.
Then for each $1 \le m \le n$ the random partition $\pi_m$ is regenerative with respect to deletion of some part $X_m$,
whose distribution $q(m, \cdot)$ is that of $H_m$ given $H_m >0 $, where $H_m$ is the number of balls
in the sample of size $m$ which fall in some particular box containing $X_n$ balls 
in $\pi_n$.
\end{theorem}

\noindent
The main point of this theorem is its implication that if
$\pi_n$ is regenerative with respect to deletion of $X_n$ according to
some deletion kernel $d(\lambda,\cdot)$, which might be defined in 
the first instance only for partitions $\lambda$ of $n$, then there is 
for each $1 \le m \le n$ an essentially unique way to construct 
$d(\lambda,\cdot)$ for partitions $\lambda$ of $m$,
so that formula \re{d-eq} holds also for $m$ instead of $n$. Iterated deletion of parts of
$\pi_n$ according to this extended deletion kernel puts the parts of $\pi_n$ in 
a particular random order, call it the {\em order of deletion according to $d$}.
This defines a random {\em composition of} $n$, that is a sequence of strictly positive integer random variables (of random length) with sum $n$. We may represent 
such a random composition of $n$ as an infinite sequence of random variables, 
by padding with zeros.
The various distributions involved in this representation of $\pi_n$
are spelled out in the following corollary, which follows easily from the
theorem.

\begin{corollary}
\label{crlq}
In the setting of the preceding theorem,

\begin{itemize}
\item [{\em (i)}]
for each $1 \le m \le n$ the distribution $q(m,\cdot)$ of $H_m$ is derived from $q(n,\cdot)$
by the formula
\begin{equation}\label{hypgeom}
q(m,k) = { q_0(m,k) \over 1 - q_0(m,0) }~~~~~~~(1 \le k \le m )
\end{equation}
where
$$
q_0(m,k) := \sum_{x = 1}^n q(n,x) { { n - x \choose m - k } { x \choose k } \over {n \choose m } }\,
~~~~~~~(0 \le k \le m ).
$$

\item [{\em (ii)}]
Let
$X_{n,1}, X_{n,2}, \ldots$ be a sequence of non-negative integer
valued random variables such that $X_{n,1}$ has distribution $q(n,\cdot)$,
and for $j \ge 1$
\eq
\label{xnj}
\prob( X_{n,j +1 } = \cdot \giv X_{n,1} + \cdots + X_{n,j} = r ) = q( n - r , \cdot )
\en
with $X_{n,j +1 } = 0$ if $X_{n,1} + \cdots + X_{n,j} = n$, so
$\XX_n := (X_{n,1}, X_{n,2}, \ldots)$ is a random composition of $n$ with
\begin{equation}\label{comp-prob}
\prob(\XX_n = \lambda)=\prod_{j=1}^\ell q( \lambda_{j}+\cdots+\lambda_{\ell},\lambda_{j})
\end{equation}
for each composition $\lambda$ 
of $n$ with $\ell$ parts of sizes $\lambda_1 ,\lambda_2 ,\ldots, \lambda_\ell$.
Then $(\pi_n,X_n)$ with the joint distribution described by {\em Theorem \ref{thm1}}
can be constructed 
as follows: let $X_n = X_{n,1}$ and define $\pi_n$ by
ranking $\XX_n$.

\item [{\em (iii)}]
For each $1 \le m \le n$ the distribution of $\pi_m$ is given by the formula
\begin{equation}\label{sumperm}
\prob(\pi_m = \lambda)=\sum_{\sigma}\prod_{j=1}^\ell q( \lambda_{\sigma(j)}+\cdots+\lambda_{\sigma(\ell)},\lambda_{\sigma(j)})
\end{equation}
where $\lambda$ is a partition of $m$ into $\ell$ parts of sizes 
$\lambda_1 \ge \lambda_2 \ge \cdots \ge \lambda_\ell >0$,
and the summation extends over all $m!/\prod a_j(\lambda)!$ distinct permutations $\sigma$ of the $\ell$
parts of $\lambda$, with $a_j(\lambda)$ being the number of parts of $\lambda$ of size $j$.
\item [{\em (iv)}]
Let $d(\lambda,x)$ for partitions $\lambda$ of $m \le n$ and $x$ a part of
$\lambda$ be derived from $q$ and $p$ via formula {\em \re{d-eq}},
and let $\XX_n$ be the random composition of $n$ defined by the parts of 
$\pi_n$ in order of deletion according to $d$. Then $\XX_n$ has the 
distribution described in part {\rm (ii)}.
\end{itemize}
\end{corollary}

Following \cite{gnedinp03}, we call a transition probability matrix 
$q(m,j)$ indexed by $1 \le j \le m \le n$, with $\sum_{j = 1}^m q(m,j) = 1$, 
a {\em decrement matrix}. 
A {\em random composition of $n$ generated by $q$} is 
a sequence of random variables $\XX_n := (X_{n,1}, X_{n,2}, \ldots)$ 
with distribution defined as in part (ii) of the previous corollary. 
Hoppe \cite{ho87} called this scheme for generating a random composition of
$n$ a \emph{discrete residual allocation model}.

Suppose now that $\XX_n$ is the sequence of sizes of classes in a 
random ordered partition $\OPi_n$ of the set $[n]$, 
meaning a sequence of disjoint non-empty sets whose union is $[n]$, 
and that conditionally given $\XX_n$ all possible choices of $\OPi_n$ 
are equally likely. 
Let $\XX_m$ be the sequence of sizes of classes of the ordered partition of $[m]$ defined
by restriction of $\OPi_n$ to $[m]$.
Then the $\XX_m$ is said to be {\em derived from $\XX_n$ by sampling}, and the
sequence of distributions of $\XX_m$ is called {\em sampling consistent.}
A {\em composition structure} is a sampling consistent 
sequence of distributions of compositions $\XX_n$ of $n$ for $n = 1,2, \ldots$. 

\begin{definition}
\label{defregcomp}
{\em Following \cite{gnedinp03},
we call a random composition $\XX_n = (X_{n,1}, X_{n,2}, \ldots)$ of $n$ {\em regenerative}, 
if for each $1 \le x < n$,
conditionally given that $X_{n,1} = x$ the remaining composition 
$(X_{n,2}, \ldots)$ of $n-x$ is distributed according to the unconditional distribution of $\XX_{n-x}$ derived from
$\XX_n$ by random sampling.
Call a composition structure $(\XX_n)$ {\em regenerative} if $\XX_n$ is regenerative for each $n = 1,2, \ldots$.
}
\end{definition}
Note the close parallel between this definition of regenerative compositions and Definition \ref{defreg}
of regenerative partitions. The regenerative property of a random partition is more subtle, because it 
involves random selection of some part to delete, and this selection process is allowed to be as general
as possible, while for random compositions it is simply the first part that is deleted.
The relation between the two concepts is provided by the following
further corollary of Theorem \ref{thm1}:

\begin{corollary}
\label{cordelorder}
If the parts of a regenerative partition $\pi_n$ of $n$ are put in 
deletion order to define a random composition of $\XX_n$ of $n$, as in
part {\em (iv)} of the previous corollary, then $\XX_n$ is a 
regenerative composition of $n$. 
\end{corollary}

This reduces the study of regenerative partitions 
to that of regenerative compositions, for which a rather complete theory has 
already been presented in \cite{gnedinp03}. 
In particular, the basic results of \cite{gnedinp03}, recalled here in Section
\ref{paintbox}, provide an explicit
paintbox representation of regenerative partition structures, along
with an integral representation of corresponding decrement matrices $q$.
See also Section \ref{corols} for some variants of
Corollary \ref{cordelorder}.

\section{Proof of Theorem \ref{thm2}}
\label{proofdk}
This is an extension of the argument of Kingman \cite{ki78a} in the
case $\tau = 0$.
Recall first that when partitions $\lambda$
are encoded by their multiplicities, $a_r = a_r(\lambda)$ for
$r = 1,2, \ldots$, the sampling consistency condition on a partition probability function $p$ 
is expressed by the formula
\eq
\label{pconsist}
p(a_1,a_2, \ldots) = p(a_1 +1,a_2, \ldots) { a_1 +1 \over n + 1 }
+ \sum_{r>1} p ( \ldots, a_{r-1} - 1, a_r + 1 , \ldots ) { r ( a_r + 1 ) \over n +1 }
\en
where $p$ is assumed to vanish except when its arguments are non-negative integers,
and $n = \sum_r r a_r$.
\par Assuming that $p$ is a regenerative with respect to $d$,
iterating (\ref{d-eq}) we have for parts $r,s\in\lambda$, 
\eq
\label{permp}
p(\lambda)={q(n,r)\over d(\lambda,r)}\,{q(n-r,s)\over d(\lambda - \{r\},s)}\,p(\lambda - \{r,s\}),
\en
which can clearly be expanded further.
Since this expression is invariant under permutations of the parts,
interchanging $r$ and $s$ we get
$${q(n,r)\over d(\lambda,r)}\,{q(n-r,s)\over d(\lambda - \{r\},s)}=
{q(n,s)\over d(\lambda,s)}{q(n-s,r)\over d(\lambda - \{s\}, r )}.$$
\par Assume now that $d$ is given by (\ref{dk}).
Introducing 
$$b(n,r):={q(n,r) n\over (n-r)\tau +r(1-\tau)}$$
formula \re{permp}
yields $b(n,r)b(n-r,s)=b(n,s)b(n-s,r)$. Taking $s=1$ and 
abbreviating $f(n):=b(n,1)$ we obtain $b(n,r)/b(n-1,r)=f(n)/f(n-r)$, thus
$$b(n,r)=f(n-r+1)\cdots f(n-1)f(n)g(r)\,, {~~\rm for ~~} g(r):={b(r,r)\over f(1)\cdots f(r)}\,.$$
The full expansion of $p$ now reads
$$p(\lambda)= \prod_{i=0}^{\ell-1}(1-\tau+i\,\tau)\prod_{k=1}^n f(k) \prod_r {g(r)^{a_r}\over a_r!}$$
where $a_r$ is the number of parts of $\lambda$ of size $r$, with $\Sigma\, a_r=\ell$ and $\Sigma \, r a_r = n$.
By homogeneity we can choose the normalisation $g(1)=1$.
Assuming that $p$ is a partition structure, substituting into 
(\ref{pconsist}) 
and cancelling common terms gives
$${n+1\over f(n+1)}= (1-\tau+\ell\,\tau)+\sum_{r>1} r\,a_{r-1} {g(r)\over g(r-1)}.$$ 
Now defining $h(r)$ by the substitution 
$${g(r)\over g(r-1)}=-{\tau\over r}+{r-1\over r}\,h(r)$$
we obtain
$$ {n+1\over f(n+1)}=1-\tau+\sum_{r>1}(r-1)a_{r-1}h(r)$$
which must hold for arbitrary partitions, hence $h(r)=\gamma$ for some constant. Therefore
$$f(n)={n\over 1-\tau+(n-1)\gamma}\,\,,\qquad g(r)={(\gamma-\tau)_{r-1\,\uparrow\gamma}\over r!}.$$
\par It follows that 
$$p(\lambda)={n!\,(1-\tau)_{\ell\,\uparrow\tau}\over (1-\tau)_{n\uparrow\gamma}}
\prod_r\left({(\gamma-\tau)_{r-1\,\uparrow\gamma}\over
r!}\right)^{a_r} {1\over a_r!}$$
which is positive for all $\lambda$ iff $\gamma>\tau$. 
The substitution 
$$ \alpha={\tau\over\gamma}\,\,,\qquad \theta={1-\tau\over\gamma}\,$$
reduces this expression to the two-parameter formula \re{param2},
and Theorem \ref{thm2} follows.

\section{Fragmented permutations}
\label{frag}

We use the term {\em fragmented permutation of $[n]$}
for a pair 
$\gamma = (\sigma, \lambda) \in \SS_n \times C_n$, where $\SS_n$ is 
the set of all permutations of $[n]$, and
$C_n$ is the set of all compositions of $n$.
We interpret a fragmented permutation $\gamma$
as a way to first arrange $n$ balls labeled by $[n]$ in a sequence,
then fragment this sequence into some number of {\em boxes}.
We may represent a fragmented permutation
in an obvious way, e.g.
$$
\gamma = 2,3,9 \boxbar 1,8 \boxbar 6,7,5 \boxbar 4
$$
describes the configuration with balls 2, 3 and 9 in that order in the
first box, balls 1 and 8 in that order in the second box, and so on,
that is $\gamma = (\sigma, \lambda)$ for
$\sigma = (2,3,9,1,8,6,7,5,4)$ and $\lambda = (3,2,3,1)$.

We now define a transition probability matrix on the set of all
fragmented permutations of $[n]$.
We assume that some probability distribution $q(n, \cdot)$ is
specified on $[n]$. Given some initial fragmented permutation
$\gamma$, 
\begin{itemize}
\item
let $X_n$ be a random variable with distribution $q(n,\cdot)$, meaning
$$
\prob (X_n = x ) = q(n,x), ~~~ 1 \le x \le n ;
$$
\item
given $X_n=x$, pick a sequence of $x$ different balls uniformly at random
from the
$$
n (n-1) \cdots (n - x + 1)
$$
possible sequences;
\item
remove these $x$ balls from their boxes and put them in order in 
a new box to the left of the remaining $n-x$ balls in boxes.
\end{itemize}
To illustrate for $n = 9$, if the initial fragmented permutation
is $\gamma = 2,3,9 \boxbar 1,8 \boxbar 6,7,5 \boxbar 4$ as above,
$X_9 = 4$ and the sequence of balls chosen is $(7,4,8,1)$, then the new
fragmented permutation is
$$
7,4,8,1 \boxbar 2,3,9 \boxbar 6,5 .
$$
\begin{definition}
{\em Call the Markov chain with this transition mechanism the
{\em $q(n,\cdot)$-chain on fragmented permutations of $[n]$.}
}
\end{definition}

To prepare for the next definition,
we recall a basic method of transformation of transition probability functions.
Let $Q$ be a transition probability matrix on a finite set $S$, and let $f: S \to T$ be a surjection from
$S$ onto some other finite set $T$.
Suppose that the $Q(s, \cdot)$ distribution of $f$ depends only on the value of $f(s)$, that is
\eq
\sum_{x : f(x) = t } Q(s,x) = \widehat{Q} ( f(s), t ), ~~~~(t \in T ) 
\en
for some matrix $\widehat{Q}$ on $T$.
The following consequences of this condition are elementary and well known:
\begin{itemize}
\item
if $(Y_n, n = 0,1,2, \ldots)$ is a Markov chain with transition matrix $Q$ and starting state $x_0$, then
$(f(Y_n), n = 0,1,2, \ldots)$ is a Markov chain with transition matrix $\widehat{Q}$ and starting state 
$f(x_0)$;
\item
if $Q$ has a unique invariant probability measure $\pi$, then
$\widehat{Q}$ has unique invariant probability measure $\widehat{\pi}$ which is
the $\pi$ distribution of $f$.
\end{itemize}
To decribe this situation, we may say that $\widehat{Q}$ is the {\em push-forward of $Q$ by $f$.}

\begin{definition}
{\em The {\em $q(n,\cdot)$-chain on permutations of $[n]$}
is the $q(n,\cdot)$-chain on fragmented permutations of $[n]$
pushed forward by projection from $(\sigma,\lambda)$ to $\sigma$.
Similarly, pushing forward from $(\sigma,\lambda)$ to $\lambda$
defines the
{\em $q(n,\cdot)$-chain on compositions of $n$} and pushing forward further
from compositions to partitions, by ranking, 
defines the {\em $q(n,\cdot)$-chain on partitions of $n$}.
}
\end{definition}

The $q(n,\cdot)$-chain on permutations of $[n]$ is the inverse of the
top-to-random shuffle studied by Diaconis, Fill and Pitman \cite{dfp92}. 
Keeping track of packets of cards in this shuffle leads naturally
to the richer state space of fragmented permutations. 
The mechanism of the $q(n,\cdot)$-chain on compositions of $n$ is identical to that described
above for fragmented permutations, except that the labels of the balls are ignored.
The mechanism of the
$q(n,\cdot)$-chain on partitions of $n$ is obtained by further ignoring the order of
boxes in the composition. The following lemma connects these Markov chains to the
basic definitions of regenerative partitions and regenerative compositions which we made in Section \ref{intro}.

\begin{lemma}
Let $q(n,\cdot)$ be a probability distribution on $[n]$. Then
\begin{itemize}
\item[{\rm (i)}] 
a random composition $\XX_n = (X_{n,1}, X_{n,2}, \ldots )$, with 
$X_{n,1}$ distributed according to $q(n,\cdot)$, is regenerative if and only if
the distribution of $\XX_n$ is an invariant distribution for the $q(n,\cdot)$-chain on compositions of $n$,
\item[{\rm (ii)}] 
a random partition $\pi_n$ is regenerative with respect to deletion of some part $X_n$ with distribution
$q(n,\cdot)$ if and only if the distribution of $\pi_n$ is an invariant distribution
for the $q(n,\cdot)$-chain on partitions of $n$.
\end{itemize}
\end{lemma}

\proof
The proofs of the two cases are similar, so we provide details only in case (ii).
The condition that $\pi_n$ is regenerative with respect to deletion
of $X_n$ can be written as follows:
$$
\pi_n - \{ X_n \} \ed \hat{\pi}_{n - X_n}
$$
where 
\begin{itemize}
\item [(a)] 
$\ed$ denotes equality in distribution of two random elements of 
the set of partitions of $m$ for some $0 \le m \le n$, allowing
a trivial partition of $0$.
\item [(b)] 
on the left side $\pi_n - \{ X_n \}$ denotes the random partition
of $n-X_n$ derived from $\pi_n$ by deletion of the part $X_n$ of $\pi_n$,
\item [(c)] 
on the right side $(\hat{\pi}_{m}, 0 \le m \le n)$ is a sampling
consistent sequence of random partitions, independent of $X_n$, 
with $\hat{\pi}_n \ed \pi_n$.
\end{itemize}
Consider the random partition
$$
\pi_n^* := \hat{\pi}_{n -X_n} \cup \{X_n \} 
$$
obtained from $\hat{\pi}_n$ by first removing $X_n$ balls from $\hat{\pi}_n$
by random sampling, then putting all these balls in a new box.
The conditional distribution of $\pi_n^*$ given $\hat{\pi}_n$ defines
a transition probability matrix on the set of partitions of $n$, which is
the transition matrix of the $q(n,\cdot)$-chain on partitions 
of $n$. 
If $\pi_n$ is regenerative with respect to deletion of $X_n$,
then $\pi_n \ed \hat{\pi}_n \ed \pi_n^*$.
That is to say, the distribution of $\pi_n$ is an invariant probability 
measure for the $q(n,\cdot)$-chain on partitions of $n$. 
Conversely, if the distribution of $\hat{\pi}_n$ is
invariant for the $q(n,\cdot)$-chain on partitions of $n$,
so $\hat{\pi}_n \ed \pi_n^*$,
we can set $\pi_n:= \pi_n^*$, and then by construction
$$
\pi_n - \{X_n\} = \pi_n ^* - \{X_n\} = \hat{\pi}_{n -X_n} .
$$
So $\pi_n$ is regenerative with respect to deletion of $X_n$ with 
distribution $q(n,\cdot)$.
\endpf

Theorem \ref{thm1} and its corollaries now follow easily from the previous
lemma and the following lemma:

\begin{lemma}
\label{l1}
For each probability distribution $q(n,\cdot)$ on $[n]$, the
$q(n,\cdot)$-chain on fragmented permutations of $[n]$
has a unique stationary distribution. Under this distribution, 
\begin{itemize}
\item [{\em (i)}]
the permutation of $[n]$ has uniform distribution;
\item [{\em (ii)}]
the permutation and the composition are independent;
\item [{\em (iii)}]
the composition of $n$ is generated by $q$ according
to {\em Corollary \ref{crlq}}.
\end{itemize}
Hence, the distribution on compositions described by
{\rm Corollary \ref{crlq} (ii)} is the unique stationary probability distribution for
the $q(n,\cdot)$-chain on compositions of $n$, and distribution of $\pi_n$
described by formula {\rm \re{sumperm}} is the unique stationary probability distribution 
for the $q(n,\cdot)$-chain on partitions of $n$. 
\end{lemma}

\proof
Let $m:= \max \{x : q(n,x) > 0 \}$, and write $n = k m + r$ for positive
integers $k$ and $r$ with $1 \le r < m$. 
We argue that whatever the initial state $\gamma$, there is a strictly
positive probability that after $k+1$ steps the
$q(n,\cdot)$-chain on fragmented permutations reaches
the state
$$
1,2, \ldots, m \boxbar \cdots \boxbar (k-1) m + 1, \ldots, k m \boxbar km+1, \ldots, n 
$$
in which the permutation is the identity and the composition
is $(m, \cdots, m, r)$.
To see this, note that the transition mechanism ensures that
after one step from $\gamma$ it is possible to reach a state of the form
$$
n - m +1, n -m +2, \ldots, n \boxbar \ldots \ldots
$$
for some $\ldots \ldots$ determined by the initial configuration $\gamma$.
Then after two steps it is possible to reach a state of the form
$$
(k-1) m + 1, \ldots, k m \boxbar km+1, \ldots, n \boxbar \ldots \ldots
$$
for some $\ldots \ldots$, and so on. Since there is a state which can be 
reached eventually no matter what the initial state, it follows from
the elementary theory of Markov chains that a stationary distribution
exists and is unique.

Let $P_{q(n,\cdot)}$ denote the push-forward of this stationary distribution 
to compositions of $n$, that is the stationary distribution of the
$q(n,\cdot)$-chain on compositions of $n$.
By definition of the $q(n,\cdot)$-chain on fragmented permutations,
\begin{itemize}
\item
under $P_{q(n,\cdot)}$ the number of balls in the first box has distribution
$q(n,\cdot)$;
\item
under $P_{q(n,\cdot)}$, for each $1 \le m < n$, given that the first box 
contains $n - m$ balls, the remaining composition of $m$ has the
distribution on compositions of $m$ derived from $P_{q(n,\cdot)}$ by
taking a random sample of $m$ out of the $n$ balls, to be denoted
$( P_{q(n,\cdot)} )^{n \rightarrow m }$.
\end{itemize}
To complete the proof of the lemma, it just remains to check the key identity
\begin{equation}
\label{consist}
( P_{q(n,\cdot)} )^{n \rightarrow m } = P_{q(m,\cdot)}
\end{equation}
where $q(m,\cdot)$ is derived from $q(n,\cdot)$ by formula \re{hypgeom}.
Due to independence of the composition and the permutation,
a composition of $m$ with distribution
$( P_{q(n,\cdot)} )^{n \rightarrow m }$ is obtained from the stationary distribution
of the $q(n,\cdot)$-chain on fragmented permutations
by ignoring balls $m+1, \ldots, n$, and considering the 
composition of $m$ induced by the balls labeled by $[m]$.
But when the fragmented permutation evolves
according to the $q(n,\cdot)$-chain, it is clear that 
at each step, no matter what the initial state,
for $0 \le s \le m$, the probability that
exactly $s$ of the $m$ balls get moved to the left is
$q_0(m,s)$ as in \re{hypgeom}.
Let $q(m, \cdot)$ be as in formula \re{hypgeom}.
Since the probability of moving at least one of the first $m$ balls 
is $1 - q_0(m,0)$,
no matter what the initial state, the
$q(n,\cdot)$-chain on fragmented permutations of $[n]$ pushes forward to
a Markov chain on fragmented permutations of $[m]$. The transition
matrix of this chain is a mixture of the identity matrix and the 
matrix of $q(m, \cdot)$-chain on fragmented permutations of $[m]$, 
with weights $q_0(m,0)$ and $1 - q_0(m,0)$
respectively. 
Hence the equilibrium distribution of this
chain with state space fragmented permutations of $[m]$ is identical to the
equilibrium distribution of the $q(m, \cdot)$-chain on fragmented permutations of $[m]$, 
whose projection onto compositions of $[m]$ is $P_{q(m,\cdot)}$.
This proves \re{consist}.
\endpf

\section{Some corollaries}
\label{corols}

This section spells out some further corollaries of Theorem \ref{thm1}.

\begin{corollary}
\label{crlqq}
The distribution of every regenerative partition $\pi_n$ of $n$ is obtained by ranking the 
components of some regenerative composition $\XX_n$ of $n$, whose distribution is uniquely
determined by that of $\pi_n$. Then $\pi_n$ is regenerative with respect to deletion of 
$X_n$ distributed like the first component of $\XX_n$. 
This correspondence, made precise by the formulae of {\em Corollary \ref{crlq}}, establishes bijections between
the following sets of probability distributions:
\begin{enumerate}
\item[{\em (i)}] 
probability distributions $q(n,\cdot)$ on $[n]$

\item[{\em (ii)}] 
distributions of regenerative compostions $\XX_n$ of $n$

\item[{\em (iii)}] 
distributions of regenerative partitions $\pi_n$ on $n$.
\end{enumerate}
\end{corollary}

\noindent
An explicit link from (iii) to (i) is provided by a recursive formula 
\cite[Equation (34)]{gnedinp03} expressing $q(n,\cdot)$ via 
the probabilities of one-part partitions $(p(j),\,j=1,\ldots,n)$.
There is also an explicit formula expressing $q(n,\cdot)$ as a rational function of
probabilities $(p(1^m),\,m=1,\ldots,n)$ where $1^m\vdash m$ is the partition with only singleton parts.

\par As a variant of the above corollary, we record also:
\begin{corollary}
Given a decrement matrix $q = (q(m,j), 1 \le j \le m \le n)$ for some fixed $n$,
for $1 \le m \le n$ let $\XX_m$ be the random composition of $m$ generated by $q$, 
and $\pi_m$ the random partition of $m$ obtained by ranking $\XX_m$. 
The following conditions are equivalent:
\begin{enumerate}
\item[{\em (i)}] 
the entire decrement matrix $q$ is determined by $q(n, \cdot)$ according to formula \rem{hypgeom};
\item[{\em (ii)}] 
the sequence of compositions $(\XX_m, 1 \le m \le n)$ derived from $q$ is sampling consistent;
\item[{\rm (iii)}] 
the sequence of partitions $(\pi_m, 1 \le m \le n)$ derived from $q$ as in \rem{sumperm} is sampling consistent.
\end{enumerate}
\end{corollary}
The equivalence of (i) and (ii) can also be read from \cite{gnedinp03}, where formula
\re{hypgeom} was given only for $m = n-1$ as a means of recursively computing $q(m,\cdot)$ from
$q(n,\cdot)$ for $m < n$.
That (ii) implies (iii) is obvious. That (iii) implies (ii) is not obvious,
but this is an immediate consequence of Theorem \ref{thm1} and Corollary \ref{crlq}, 
because (iii) means that $\pi_n$ is regenerative with respect to deletion of
the first term of $\XX_n$.

\begin{corollary}
\label{crl5}
The following two conditions on a random composition 
$\XX_n = (X_{n,1}, X_{n,2}, \ldots)$ of $n$ are equivalent:
\begin{enumerate}
\item[{\em (i)}] 
$\XX_n$ is regenerative;
\item[{\em (ii)}] 
for each $1 \le j \le k < n$, conditionally given 
$X_{n,1}, \ldots, X_{n,j}$
with 
$X_{n,1} + \cdots + X_{n,j} = k$,
the partition of $n-k$ obtained by ranking
$X_{n,j+1}, X_{n,j+2}, \ldots$ has the same distribution as
$\pi_{n-k}$, the partition of $n-k$ obtained by sampling from $\XX_n$.
\end{enumerate}
\end{corollary}
\proof
That (i) implies (ii) is obvious. Conversely, condition (ii) for
$j = 1$ states that the partition $\pi_n$ derived from $\XX_n$
is regenerative with respect to deletion of $X_{n,1}$. Let 
$q(n,\cdot)$ be the distribution of $X_{n,1}$, and let $d$
denote the corresponding deletion kernel, extended to partitions
of $m$ for $1 \le m \le n$ in accordance with Theorem \ref{thm1}.
Condition (ii) for $j = 2$ implies that for each $i$ such that
$q(n,i) > 0$, 
the partition $\pi_{n-i}$ is regenerative with respect to deletion of a 
part whose unconditional distribution equals the conditional distribution 
of $X_{n,2}$ given $X_{n,1} = i$. According to the uniqueness
statement of Theorem \ref{thm1}, this distribution must be the
distribution $q(n-i, \cdot)$ determined by $q(n,\cdot)$
via \re{hypgeom}. So the joint law of $(X_{n,1}, X_{n,2})$ is
identical to that described in Corollary \ref{crlq} (ii).
Continuing in this way, it is clear that the distribution of the
entire sequence $\XX_n$ is that described in Corollary \ref{crlq} (ii).
The conclusion now follows from Corollary \ref{cordelorder}.
\endpf

\section{Paintbox representations}
\label{paintbox}

Gnedin's paintbox representation of composition structures \cite{gnedin97} uses a random closed 
set ${\cal R}\subset [0,1]$ to separate points of a uniform sample into clusters. 
Given $\cal R$, define an interval partition of $[0,1]$ comprised of
{\em gaps}, that is open interval components of $[0,1]\setminus {\cal R}$,
and of individual points of ${\cal R}$.
A random ordered partition of $[n]$ is then constructed from ${\cal R}$ and 
independent uniform sample points $U_1,\ldots,U_n$ 
by grouping the indices of sample points which fall in the same gap, and letting the points which hit 
$\cal R$ to be singletons.
A random composition $\XX_n$ of $n$ is then constructed as the sequence of block sizes in this
partition of $[n]$, ordering the blocks from left to right, according to the location of the
correspoding sample points in $[0,1]$.
Gnedin showed that every composition structure $(\XX_n)$ can be so represented.
As in Kingman's representation of partition structures, $\cal R$ can be interpreted as an asymptotic shape of 
$\XX_n$, provided $\XX_n$ is properly encoded as an element of the metric space of closed subsets of $[0,1]$ 
with the Hausdorff distance function.
\par
According to the main result of \cite{gnedinp03}, each regenerative composition structure $(\XX_n)$ 
is associated in this way with an $\cal R$ which is {\em multiplicatively regenerative} in the
following sense: for $t \in [0,1]$ let $D_t:= \inf ([t,1] \cap {\cal R })$, and given that
$D_t < 1$ let 
$$
{\cal R}^{[D_t,1]} := \{ (z - D_t )/(1 - D_t), z \in {\cal R}\cap [D_t,1]\}
$$
which is the restriction of ${\cal R}$ to $[D_t,1]$ scaled back to $[0,1]$;
then
\eq
\label{mregen}
( {\cal R}^{[D_t,1]} \giv D_t \mbox{ with } D_t < 1\,,{\cal R}\cap [0,D_t]) \ed {\cal R }
\en
meaning that the conditional distribution of ${\cal R}^{[D_t,1]}$, given $D_t$ with $D_t < 1$ and 
given ${\cal R}\cap [0,D_t]$,
is identical to the unconditional distribution of ${\cal R }$.
This condition holds if and only if $- \log ( 1 - {\cal R} )$ is a regenerative
subset of $[0,\infty[$ in the usual (additive) sense. So by a result of Maisonneuve, 
the most general multiplicatively regenerative subset ${\cal R}$ can be constructed as the
closure of $\{ 1 - \exp(-S_t) , t \ge 0 \}$ for some subordinator $(S_t, t \ge 0 )$.
Thus regenerative composition structures are parameterised by a pair $(\nutil,\drift)$ where
$\nutil$ is a measure on $\,]0,1]$ with finite first moment and 
$\drift\geq 0$. 
The measure $\nutil({\rm d}u)$ is the image of the L{\'e}vy measure $\nu({\rm d}s)$ of the
subordinator via the transformation from $s$ to $1 - \exp(-s)$.
and $\drift$ is the drift parameter of the subordinator. So the Laplace
exponent of the subordinator, evaluated at a positive integer $n$,
is
$$
\Phi(n)=n\drift+\int_{]0,1]}(1-(1-x)^n)\nutil({\rm d}x)\,,\qquad
$$
The decrement matrix of the regenerative composition is then
$$
q(n,r)={\Phi(n,r)\over \Phi(n)}\,,\qquad 1\leq r\leq n\,,~n=1,2,\ldots
$$
where
$$
\Phi(n,r)=n\drift\,1(r=1)+ {n\choose r}\int_{]0,1]} x^r(1-x)^{n-r}\nutil({\rm d}x)\,.
$$ 
Uniqueness of the parametrisation is achieved by a normalisation condition, e.g. $\Phi(1)=1$. 
\vskip0.5cm
\noindent

The partition structure derived by sampling from a random closed subset ${\cal R}$ of $[0,1]$
depends only on the distribution of the sequence of {\em ranked lengths induced by} $\cal R$
$$
V({\cal R}) := (V_1({\cal R}), V_2({\cal R}), \ldots )
$$
where $V_i({\cal R})$ is the length of the $i$th longest interval component of $[0,1]\setminus \cal R $.
Our consideration of regenerative partition structures suggests the following definition.
Call $\cal R$ {\em weakly multiplicatively regenerative} if
for each $t \in [0,1]$ 
\eq
\label{wmregen}
( V( {\cal R}^{[D_t,1]} ) \giv D_t \mbox{ with } D_t < 1\,,{\cal R}\cap [0,D_t]) \ed V({\cal R })
\en
meaning that the conditional distribution of relative ranked lengths induced by 
${\cal R} \cap {[D_t,1]}$, given $D_t$ with $D_t < 1$ and given the restricted set 
${\cal R}\cap [0,D_t]$, 
is identical to the unconditional
distribution of ranked lengths induced by $\cal R$.
From Theorem \ref{thm1} we easily deduce:

\begin{corollary}
A random closed subset $\cal R$ of $[0,1]$ is weakly multiplicatively regenerative if and only if
it is multiplicatively regenerative.
\end{corollary}
\proof
The ``if'' part is obvious by measurability of the map from $\cal R$ to $V( { \cal R } )$.
To argue the converse, suppose that $\cal R$ is 
weakly multiplicatively regenerative. Without loss of generality, it can be supposed that
$\cal R$ is defined on the same probability space as a sequence of independent uniform
$[0,1]$ variables $U_i$ for $i = 1,2, \ldots$.
Let $\XX_n = (X_{n,1}, X_{n,2}, \ldots )$ for $n = 1,2, \ldots$ be the sequence
of compositions of $n$ derived from $\cal R$ by sampling with these independent uniform variables,
and let $\pi_n$ be the partition of $n$ defined by ranking $\XX_n$.
By consideration of \re{wmregen} with $t$ replaced by $U_{n,1}$,
where $U_{n,k}$ is the $k$th order statistic of $U_1,\ldots,U_n$,
it is easily argued that $\pi_n$ is regenerative with respect to deletion of $X_{n,1}$.
If $X_{n,1}<n$, we apply (\ref{wmregen}) at
$t=U_{n,X_{n,1}+1}$ and repeat the argument to show that 
$\pi_n-\{X_{n,1},X_{n_2}\}$ is a distributional copy of $\pi_{n-m}$ given 
$X_{n,1}$ and $X_{n,2}$ with $X_{n,1}+X_{n,2}=m$.
Iterating further we see that $\XX_n$ satisfies condition (ii) of Corollary
\ref{crl5}. Hence $\XX_n$ is regenerative. Thus $(\XX_n)$ defines
a regenerative composition structure, and it follows the main result 
of \cite{gnedinp03} that the set ${\cal R}$ is regenerative. 
\endpf

\section{Uniqueness and positivity}

For each $n$ equation (\ref{d-eq}) is the basic algebraic relation between the functions $p(\lambda)$, $d(\lambda)$ 
for $\lambda\vdash n$ and $q(n,\cdot)$. According to Corollary \ref{crlq} $q(n,\cdot)$ uniquely determines $p$,
for each $n$, but $d$ satisfying (\ref{d-eq}) need not be unique because $p(\lambda)$ may assume zero value for some $\lambda$, as illustrated by
the following example. 
\vskip0.5cm
\noindent
{\bf Example} (Regenerative hook partition structures.)
A partition $\lambda$ of $n$ is called a {\em hook} if 
it has at most one part larger than $1$.
The only regenerative partition structures $p$ such that
$\pi_n$ is a hook with probability one for every $n$ are those which
can be generated by
$\nutil=\delta_1$, the Dirac mass at $1$ and
$\drift\geq 0$, including the trivial boundary case with $\drift=\infty$. 
Then for $n>1$ 
$$q(n,n)={1\over 1+n\drift}\,,\qquad q(n,1)={n\drift\over 1+n\drift}.$$
This implies that the
associated {\it hook composition structure} gives positive probability only to compositions of the type
$(n)$ or $(1^m,n-m)$, $1\leq m\leq n$, for every $n$. 
For these hook composition structures 
the deletion kernel is an arbitrary kernel with the property
\begin{equation}\label{d-hook}
~d(\lambda,\,1)=1\,~{\rm if~}1\in\lambda\,.
\end{equation}
This property is characteristic, that is each partition structure regenerative according to such deletion kernel $d$ is derived from a hook composition structure.
Indeed, if a composition structure compatible with such $d$ gave positive probability
to a composition $(x,y,\ldots)$ with $x>1$ then, by sampling consistency, the
composition $(2,1)$ would also have positive probability, contradicting $q(3,2)=0$.
\vskip0.5cm
\par 
The next lemma shows that only in the hook case can there be any ambiguity about the deletion kernel 
generating some partition structure 
\begin{lemma}\label{posit}
Let $(\pi_n)$ be a regenerative partition structure. Then 
\begin{enumerate}
\item[{\rm (i)}] either $(\pi_n)$ is a hook partition structure,
\item[{\rm (ii)}] or $p>0$, $d>0$ and $q>0$ for all admissible values of arguments. 
\end{enumerate}
Moreover, {\rm (i)} holds iff $p(2,2)=0$, while {\rm (ii)} holds iff $p(2,2)>0$.
\end{lemma}
{\it Proof.} If $V_2$ in Kingman's representation is strictly positive with nonzero probability, then
$p(n,n)>0$ for all $n$. By virtue of $p(n,n)=q(2n,n)p(n)$ follows $q(2n,n)>0$, but then by Corollary \ref{crlq}(i)
$q(n,\cdot)>0$ and this implies $p>0$ and $d>0$.
In this case $p(2,2)>0$.
Alternatively, if ${\mathbb P}(V_2=0)=1$ then $p(2,2)=0$ and we are in the hook case.
\endpf
\noindent
Thus, if $d$ is not a kernel satisfying (\ref{d-hook}), then
$d(\lambda,x)=0$ for some partition $\lambda$ and some
part $x\in \lambda$ implies that a
partition structure regenerative according to $d$ is trivial, 
({\it i.e.} is either the pure-singleton partition with $p(1^n)\equiv 1$, or 
the one-block partition with $p(n)\equiv 1$).

\par The positivity condition in lemma rules out nontrivial partition structures, which have
an absolute bound on the number of parts for all $n$.
For example, none of the members of the 
two-parameter family of partition structures (\ref{param2}) 
is regenerative for
$\alpha<0$ and $-\theta/\alpha =k\in \{2,3,\ldots\}$
because each $\pi_n$ has at most $k$ parts.

\begin{corollary} If a partition structure is regenerative and satisfies $p(2,2)>0$ then 
$q$ uniquely determines $p$ and $d$, and $p$ uniquely determines $q$ and $d$. Thus if a regenerative partition 
structure is not hook the corresponding deletion kernel is unique. 
\end{corollary}


\par Checking if 
a given partition structure is regenerative according to an
unknown deletion kernel can be done by first computing $q$, by some algebraic manipulations,
from a given partition probability function $p$, then computing a partition probability function $p_*$ for
the regenerative partition structure related to this $q$,
and finally checking if $p=p_*$. 
When this method is applied to a partition from the two-parameter family
with $0<\alpha<1,\,-\alpha<\theta<0$, the resulting $q$ recorded in \cite[Equation (39)]{gnedinp03}
is not everywhere positive, hence the partition structures with such parameters are not regenerative.

\section{Generalisations and related work}
\label{related}

Given some deletion kernel $d$, 
it is of interest to consider
pairs of partition structures $(p_0,p_1)$ such that 
$d$ {\em reduces $p_0$ to $p_1$}, meaning 
that the following extension of formula \re{d-eq} holds:
\begin{equation}\label{d-eqz}
p_0(\lambda) d(\lambda,x) = q(n,x)p_1(\lambda - \{x\})\,,\qquad x\in \lambda
\qquad ( \lambda \vdash n )
\end{equation}
where 
\eq
q(n,x):=\sum_{\{\lambda\vdash n\,:\,x\in \lambda\}} d(\lambda,x)p_0(\lambda)
\qquad ( 1\leq x\leq n )
\en
is the unconditional probability that the deletion rule removes a part of size $x$ from $\pi_n$ distributed according to $p_0$.
\par 
Pitman \cite{jp.isbp} showed that if $p_0,p_1,p_2, \ldots$ is a sequence of
partition structures such that size-biased deletion reduces $p_i$ to 
$p_{i+1}$ for each $i \ge 0$, and $p_0$ can be represented in terms of
random sampling from $(V_1,V_2, \ldots)$ with $V_i>0$ for each $i$, then
$p_0$ is an $(\alpha,\theta)$ partition as in \re{param2}
for some $0 \le \alpha < 1$ and $\theta > - \alpha$, in which case
$p_j$ is the $(\alpha,\theta + j \alpha)$ partition structure.
This result and Theorem \ref{thm2} are two different two-parameter generalisations of 
Kingman's characterisation of $(0,\theta)$ partition structures. 
In the result of \cite{jp.isbp}, the deletion kernel is still
defined by size-biased sampling, and repeated deletions 
generate a succession of partition structures. Whereas in Theorem \ref{thm2}
the deletion kernel is modified, and repeated deletions generate the same partition structure.

\par A class of partition structures satisfying (\ref{d-eqz}) 
is associated with Markovian compositions introduced in \cite{gpselfsim04}.
A composition structure of this type is derived from a 
set ${\cal R}$ which has a special leftmost interval $[0,X]$, and otherwise ${\cal R}\cap [X,1]$ is a scaled copy
of some other multiplicatively regenerative set ${\cal R}'$, which is independent of 
$X$ and has the property (\ref{mregen}).
For example, the members of the two-parameter family with $0<\alpha<1, ~\theta>-\alpha$ can be associated
with such Markovian composition structures. 
In general, such ${\cal R}$ can be represented as a transformed range of a 
finite-mean subordinator with arbitrary initial distribution.

\par This discussion leaves open a number of interesting questions.
One posed in \cite{jp.ew,csp}, and apparently still open, is the
problem of describing all pairs of partition structures $(p_0,p_1)$
such that $p_0$ reduces to $p_1$ by size-biased deletion.
One could ask the same question for other deletion kernels too. But
size-biased deletion is of special interest because of its 
natural interpretation in terms of Kingman's paintbox construction from ranked
random frequencies $(V_1, V_2, \ldots)$: if $X_n$ is
a sized biased pick from $\pi_n$ derived by sampling from such $(V_i)$ associated
with the partition structure $p_0$, 
then $X_n/n$ converges in distribution to $\tilde{V_1}$ which is a size-biased
pick from the limiting frequencies, meaning that
$$
\prob( \tilde{V} = V_i \giv V_1, V_2, \ldots ) = V_i
$$
$$
\prob( \tilde{V} = 0 \giv V_1, V_2, \ldots ) = 1 - \Sigma_{i} V_i
$$
where it is assumed for simplicity that the $V_i$ are almost surely distinct.
The probability distribution of $\tilde{V}$ on $[0,1]$, known as the {\em structural
distribution} encodes many important features of the partition structure $p_0$. In
particular, $p_0(n) = \ex ( \tilde{P}^{n-1} )$ for $n = 1,2, \ldots$ (see \cite{jp.isbp,csp, gpselfsim04}).
See also \cite{py95rdd, gpselfsim04} for related work.

\paragraph{Central measures}
We comment briefly on how our results link to the potential theory on graded graphs, 
as developed by the Russian school in connection with the asymptotic representation theory of the symmetric group, see
\cite{kerov96} for a survey.
\par For $\lambda\vdash n,~\mu\vdash n+1$ the number of ways to 
extend a given partition of $[n]$ with shape $\lambda$ to some partition of $[n+1]$ with shape $\mu$
does not depend on the choice of partition of $[n]$.
Denote this number $\kappa(\lambda,\mu)$. Letting ${\cal P}_n$ denote the set of partitions of $n$, consider 
the graded graph with the vertex set ${\cal P}=\bigcup_n {\cal P}_n$
and multiplicity $\kappa(\lambda,\mu)$
for the edge connecting $\lambda$ and $\mu$. 
A partition structure
is a nonnegative solution of the recursion
$$p(\lambda)=\sum_{\mu\vdash n+1} \kappa(\lambda,\mu)p(\mu)\,,~~~~\lambda\vdash n\,,~~n=1,2,\ldots$$
with $p(1)=1$.
\par 
For every $\lambda\vdash n$ a directed path in ${\cal P}$ connecting the root $(1)$ and $\lambda$
encodes a partition of $[n]$ with shape $\lambda$. Thus $p$ defines a 
a {\it central measure} on the space of infinite paths of ${\cal P}$, with the property
that all paths of length $n$, connecting $(1)$ and $\lambda\vdash n$, have the same 
probability depending only on $\lambda$, 
for every $n$ and $\lambda\in {\cal P}_n$.
Each infinite path corresponds to a partition of the set ${\mathbb N}$ and the `centrality' of a measure
just means that the random partition of ${\mathbb N}$ is exchangeable.

\par Originally, the concept of central measure emerged in connection with the Young lattice ${\cal Y}$, which 
has the same vertex set as ${\cal P}$, but a different multiplicity function (dictated by the branching
rule for characters of the symmetric group). According to a classical result of Thoma,
central measures on ${\cal Y}$ can be described in terms of two random decreasing sequences $(V_j), (W_j)$
(as to be compared with a single decreasing sequence $(V_j)$ in Kingman's representation).
See \cite{KOV} for particular examples of central measures on ${\cal Y}$, and \cite{BO} for
a general construction of central measures on graphs which generalise both ${\cal Y}$ and ${\cal P}$.
\par Although the graphs ${\cal Y}$ and ${\cal P}$ are different, there is a way to push a partition structure on
${\cal P}$ to a central 
measure on ${\cal Y}$, analogous to the transition from monomial symmetric functions to Schur symmetric functions.
It would be interesting to study the image of regenerative partition structures under this mapping 
in some detail.

\paragraph{General deletion kernels}
For partitions $\mu$ and $\lambda$, interpreted as distributions
of unlabeled balls in unlabeled boxes, write $\mu \subset \lambda$ if
$\mu$ can be obtained from $\lambda$ by deletion of one or more boxes.
Call a function $d$ on pairs of partitions $(\lambda, \,\mu)$
with $\mu \subset \lambda$,
a {\it deletion kernel} if it satisfies
$$0\leq d\leq 1\,,\qquad \sum_{\lambda^-:\,\lambda^-\subset\,\lambda} d(\lambda\,,\,\lambda^-)=1\,.$$
For each $\lambda$ the deletion kernel defines a random split in two subpartitions 
$\lambda^+ \subset \lambda$ (which remains) and $\lambda^-\subset \lambda$ (which is removed).
We can also regard $d$ as transition probability function for a Markov chain which jumps from larger partition
$\lambda$ to a smaller $\lambda^+$.

\par Given $d$, for a sequence of random partitions $(\pi_n\vdash n,~ n=1,2,\ldots)$ -- which in general need not satisfy any consistency 
condition besides that the $\pi_n$'s are defined on the same probability space --
with each $\pi_n$ we associate
a split into $\pi_n^+$ and $\pi_n^-$ (with $\pi_n^-$ removed). This suggests various concepts
of invariance, depending on how much conditioning is allowed.
Thus we could distinguish a {\it weak deletion property}
\begin{equation}\label{delweak}
{\mathbb P}(\pi_n^+=\lambda^+\,|\,\pi_n^-=\lambda^-)=
{\mathbb P}(\pi_{n-r}=\lambda^+)\,,\qquad r=|\lambda^-|\,,
\end{equation}
and a {\it strong deletion property}
$${\mathbb P}(\pi_n^+=\lambda^+\,|\,|\pi_n^-|=r)={\mathbb P}(\pi_{n-r}=\lambda^+)$$
where we write $|\lambda|$ for $\Sigma\lambda_j$.
It would be interesting to extend results of this paper to such
general deletion operations.

\section{Acknowledgement} The work on this paper started during our joint visit to the
Laboratoire de Probabilit{\'e}s, University Paris VI, in March 2004. 
We thank Jean Bertoin and Marc Yor for their hospitality and support.

\def\cprime{$'$} \def\polhk#1{\setbox0=\hbox{#1}{\ooalign{\hidewidth
\lower1.5ex\hbox{`}\hidewidth\crcr\unhbox0}}} \def\cprime{$'$}
\def\cprime{$'$} \def\cprime{$'$}
\def\polhk#1{\setbox0=\hbox{#1}{\ooalign{\hidewidth
\lower1.5ex\hbox{`}\hidewidth\crcr\unhbox0}}} \def\cprime{$'$}
\def\cprime{$'$} \def\polhk#1{\setbox0=\hbox{#1}{\ooalign{\hidewidth
\lower1.5ex\hbox{`}\hidewidth\crcr\unhbox0}}} \def\cprime{$'$}
\def\cprime{$'$} \def\cydot{\leavevmode\raise.4ex\hbox{.}} \def\cprime{$'$}
\def\cprime{$'$} \def\cprime{$'$} \def\cprime{$'$}

\end{document}